\documentclass[a4paper,12pt]{article}
\usepackage{amssymb} \usepackage{amsfonts}
\usepackage{amsmath}%\usepackage{amsart}
\usepackage[utf8x]{inputenc}
\usepackage{tikz,relsize,mathtools}
\usepackage[cmtip,arrow]{xy}
\usepackage{pb-diagram,pb-xy}
\usepackage{indentfirst}
\newcommand{\pe}{\hspace*{\fill}\vspace*{1ex}$\Box$}
\usetikzlibrary{arrows}
\begin{document}

\begin{center}
\textbf{\Large Tangent bundles of Hantzsche-Wendt manifolds}\\[0.5cm]

\emph{A. G\c{a}sior~$^*$, A. Szczepa\'nski~$^{**}$}\\[0.5cm]

$^*$  Institute of Mathematics, Maria Curie-Sk{\l}odowska University, pl.~ M.~Curie-Sk{\l}odowskiej~1, 20-031 Lublin, POLAND, anna.gasior@poczta.umcs.lublin.pl\\%[0.2cm]
$^{**}$ Institute of Mathematics, University of Gda\'nsk,\newline ul.~Wita~Stwosza~57, 80-952 Gda\'nsk, POLAND,\newline matas@univ.gda.pl
\end{center}\hskip5mm

\date{\today}
%\maketitle
\newcommand{\F}{{\mathbb F}}
\newcommand{\Z}{{\mathbb Z}}
\newcommand{\Q}{{\mathbb Q}}
\newcommand{\R}{{\mathbb R}}
\newcommand{\C}{{\mathbb C}}
\newcommand{\h}{{\mathbb H}}
\newcommand{\N}{{\mathbb N}}
\date{\today}

\begin{center}
\parbox{12.2cm}{ \footnotesize \textbf{Abstract:}
We formulate a condition for the existence of a $\operatorname{Spin}^{\C}$-structure
on an oriented flat manifold $M^n$ with $H^2(M^n,\R) = 0.$ We prove that
$M^n$ has a $\operatorname{Spin}^{\C}$-structure if and only if there exist a homomorphism
$\epsilon:\pi_1(M^n)\to \operatorname{Spin}^{\C}(n)$ such $\bar{\lambda}_n\circ\epsilon =h,$
where $h:\pi_{1}(M^n)\to \operatorname{SO}(n)$ is a holonomy homomorphism and
$\bar{\lambda}_n: \operatorname{Spin}^{\C}(n)\to \operatorname{SO}(n)$ is a standard homomorphism
defined on page 2.  As an application we shall prove that all
cyclic Hantzsche - Wendt manifolds do not have the $\operatorname{Spin}^{\C}$-structure.}
\end{center}\hskip5mm


\section{Introduction}

Let $M^n$ be a flat manifold of dimension $n.$
By definition, this is a compact connected Riemannian manifold without boundary
with sectional curvature equal to zero. From the Bieberbach theorems (\cite{Ch}, \cite{S3})
the fundamental  group $\pi_{1}(M^{n}) = \Gamma$ of $M^n$
 determines a short exact sequence:
\begin{equation}\label{ses}
0 \rightarrow \Z^{n} \rightarrow \Gamma \stackrel{h}\rightarrow
F \rightarrow 0,
\end{equation}
where $\Z^{n}$ is a torsion free abelian group of rank $n$ and $F$
is a finite group which is isomorphic to the holonomy group of
$M^{n}.$ The universal covering of $M^{n}$ is the Euclidean space
$\R^{n}$ and hence $\Gamma$ is isomorphic to a discrete cocompact
subgroup of the isometry group $\operatorname{Isom}(\R^{n}) =
\operatorname{O}(n)\ltimes\R^n = \operatorname{E}(n).$ In the above
short exact sequence $\Z^{n} \cong (\Gamma \cap \R^{n})$ and $h$ can
be considered as the projection $h: \Gamma \rightarrow F\subset
\operatorname{O}(n) \subset \operatorname{E}(n)$ on the first
component. Conversely, having a short sequence of the form
(\ref{ses}), it is known that the group $\Gamma$ is (isomorphic to)
a Bieberbach group if and only if $\Gamma$ is torsion free. \vskip
1mm By Hantzsche-Wendt manifold (for short HW-manifold) $M^n$ we
shall understand any oriented flat manifold of dimension $n$ with a
holonomy group $(\Z_2)^{n-1}.$ It is easy to see that $n$ is always
an odd number. Moreover, (see \cite{MR} and \cite{S0}) HW-manifolds
are rational homology spheres and its holonomy representation\footnote{That is a representation $\phi_{\Gamma}:F\to \operatorname{GL}(n,\Z),$
given by a formula $\phi_{\Gamma}(f)(z) = \bar{f}z\bar{f}^{-1}$,
where $\bar{f}\in\Gamma, f\in F, z\in\Z^n$ and $p(\bar{f}) = f.$} is
diagonal, \cite{RS}. Hence $\pi_{1}(M^n)$ is generated by $\beta_i =
(B_i,b_i)\in \operatorname{SO}(n)\ltimes\R^n, 1\leq i\leq n,$ where
\begin{equation}\label{genhw}
B_i = \text{diag}(-1,-1,...,-1,\underbrace{1}_i,-1,-1,...,-1)\hskip 2mm \text{and}\hskip 2mm b_i\in\{0,1/2\}^{n}.
\end{equation}
Let us recall some other properties of $M^n.$ For $n\geq 5$ the
commutator subgroup of the fundamental group is equal to the
translation subgroup ($[\Gamma,\Gamma] = \Gamma\cap\R^n$), (\cite{P}).
The number $\Phi(n)$ of affinian not equivalent HW-manifolds of
dimension $n$ growths exponentially, see \cite[Theorem 2.8]{MR} and
for $m\geq 7$ there exist many pairs of isospectral manifolds that are
not homeomorphic to each other, \cite[Corollary 3.6]{MR}. These
manifolds have an interesting connection with Fibonacci groups
\cite{S1} and the theory of quadratic forms over the field $\F_2,$
\cite{S2}. HW-manifolds have not a $\operatorname{Spin}$-structure, \cite[Example 4.6 on page 4593]{MP}. Hence tangent
bundles of HW-manifolds are non-trivial.
There are yet unknown their (co)homology groups with
coefficients in $\Z.$ Here we send readers to \cite{DP} where there are presented results for low
dimensions and an algorithm. Finally, let us mention
properties related to the theory of fixed points. HW-manifolds
satisfy the so called Anosov relation. This means that any continuous map
$f:M^n\to M^n$ satisfies $\mid L(f)\mid = N(f),$ where $L(f)$ is the Lefschetz
number of $f$ and $N(f)$ is the Nielsen number of $f,$ see
\cite{DR}.

In this note we are interested in properties of the tangent bundle
of HW-manifolds. We shall prove that they are line element
parallelizable (Proposition 1) and we shall define an infinite
family of HW-manifolds without $\operatorname{Spin}^{\C}$-structure
(Theorem~\ref{maintheorem}). However, the main result of this
article is related to the existence of $\operatorname{Spin}^{\C}$-structures on oriented flat manifolds. The group
$\operatorname{Spin}^{\C}(n)$ is given by
$\operatorname{Spin}^{\C}(n) = (\operatorname{Spin}(n)\times
S^1)/\{1,-1\}$ where $\operatorname{Spin}(n)\cap S^1 = \{1,-1\}.$
Moreover, there is a homomorphism of groups $\bar{\lambda}_n:
\operatorname{Spin}^{\C}(n)\to \operatorname{SO}(n)$ given by
$\bar{\lambda}_{n}[g,z] = \lambda_n(g),$ where $g\in
\operatorname{Spin}(n), z\in S^1$ and
$\lambda_n:\operatorname{Spin}(n)\to \operatorname{SO}(n)$ is the
universal covering. We shall prove: \vskip 1mm \noindent
{\bf Theorem~\ref{main}}
{\it  Let $M$ be an oriented flat manifold with $H^2(M,\R) =
0.$ $M$ has a $\operatorname{Spin}^{\C}$-structure if and only if there exists a
homomorphism $\epsilon:\Gamma\to \operatorname{Spin}^{\C}(n)$ such that
\begin{equation}
\bar{\lambda}_n\circ\epsilon = h.
\end{equation}}
\vskip 2mm
\noindent
As an application we prove Theorem~\ref{maintheorem}.
\vskip 1mm
\noindent
{\bf Theorem~\ref{maintheorem}}
{\it All cyclic {\rm HW}-manifolds have not the $\operatorname{Spin}^{\C}$-structure}.
\vskip 2mm
\noindent

For a description of cyclic HW-manifolds see \cite[page 11]{MR} and  Definition~\ref{cyclic}.
We conjecture that all HW-manifolds
have not the $\operatorname{Spin}^{\C}$-structure.
\vskip 2mm

\section{Hantzsche-Wendt manifolds are line element parallelizable}
We keep notation from the introduction.
For any discrete group $G,$ we have a universal principal $G$-bundle with the total space $\operatorname{E}G$ and the base space $\operatorname{B}G.$
$\operatorname{B}G$ is called the classifying space of a group $G$ and is unique up to homotopy.
In our case $\R^n$ is the total space of a principal $\Gamma$-bundle with a base space $M^n.$
Here $\operatorname{E}\Gamma = \R^n$ and $\operatorname{B}\Gamma = M^n,$ see \cite[page 369]{ATV}.
Now $G\to \operatorname{B}G$ behaves more or less like a functor, and in particular, from
the surjection $h:\Gamma\to h(\Gamma) = F$ we can construct a corresponding map
$B(h): \operatorname{B}\Gamma\to \operatorname{B}F.$ Finally, the inclusion $i_{n}:F\to \operatorname{O}(n)$ yields a map
$B(i_{n}):\operatorname{B}F\to \operatorname{B}(\operatorname{O}(n)).$ The universal $n$-dimensional vector bundle over $\operatorname{B}(\operatorname{O}(n))$ yields,
via this map a vector bundle $\eta_n$ over $\operatorname{B}F.$
\newtheorem{lem}{Lemma}
\begin{lem}\label{mainprop}
{\em (\cite[Proposition 1.1]{ATV})} $B(h)^{\ast}(\eta_n)$ is equivalent
to the tangent bundle of $M^n.$
\end{lem}
{\bf Proof:} (See \cite[page 369]{ATV}) We have the following commutative diagram
$$\begin{diagram} \node{\mathbb
R^n=\operatorname{E}\Gamma}\arrow{s}\arrow{e,t}{E(h)}\node{\operatorname{E}F}\arrow{s}\\
\node{M^n=\operatorname{B}\Gamma}\arrow{e,t}{B(h)}\node{\operatorname{B}F}
\end{diagram}$$
where $E(h)(g\cdot e) = h(g)\cdot E(h)(e)$ for all $g\in\Gamma$ and $e\in \operatorname{E}\Gamma = \R^n.$
Let the total space of $\eta_n$ be $\operatorname{E}F\times\R^n/F$ where $f\in F$
acts via $f(e,v) = (f\cdot e,f\cdot v).$
Now clearly the total space $\tau$ of the tangent bundle of $M^n = \operatorname{B}\Gamma$
can be taken to be $\R^n\times\R^n/\Gamma$
where $\Gamma$ acts via $g(v_1,v_2) = (gv_1, h(g)v_2).$
Thus we have a commutative diagram as follows:
$$\begin{diagram}
\node{\tau=\mathbb R^n\times\mathbb
R^n/\Gamma}\arrow{s}\arrow{e}\node{\operatorname{E}F\times\mathbb R^n/F}\arrow{s}\\
\node{M^n=\mathbb R^n/\Gamma}\arrow{e,t}{B(h)}\node{\operatorname{B}F}
\end{diagram}$$
where $F$ acts on $\operatorname{E}F\times\R^n$ as following $\{v_1,v_2\}\to\{E(h)(v_1),v_2\}.$ This finishes the proof.

\pe
\vskip 3mm
\newtheorem{rem}{Remark}
\begin{rem}
{\em From the above Lemma we can derive that the tangent bundle is flat in the sense of \cite[page 272]{AS}.}
\end{rem}
Let us present the main result of this section.
\newtheorem{theo}{Theorem}
\newtheorem{prop}{Proposition}
\begin{prop}\label{mtheo}
Let $M^n$ be a {\rm HW}-manifold of dimension $n.$ Then its tangent bundle
is line element parallelizable, {\rm(}is a sum of line bundles{\rm)}.
\end{prop}
{\bf Proof:}
By definition the fundamental group $\Gamma = \pi_{1}(M^n)$ is a subgroup of
$\operatorname{SO}(n)\ltimes\R^n$ and
$h(\Gamma) = (\Z_2)^{n-1}\subset \operatorname{SO}(n)$ is a group of all diagonal orthogonal matrices.
It is also an image of the holonomy representation $\phi_{\Gamma}:(\Z_2)^{n-1}\to \operatorname{SO}(n).$
Let us recall basic facts on line bundles. It is well known that the classification
space for line bundles is $\R P^{\infty},$ the infinite projective space.
Hence any line bundle $\xi:L\to M^n$ is isomorphic to $f^{\ast}(\eta_1),$
where
$$f\in [M^n,\R P^{\infty}]\simeq H^1(M^n,\Z_2)\simeq Hom(\Gamma,\Z_2)\stackrel{(\ast)}\simeq (\Z_2)^{n-1}$$
is a classification map
and $\eta_1\in H^{1}(\R P^{\infty},\Z_2) = \Z_2$ is
not a trivial element. Here $\eta_1$ represents the universal line vector bundle and the isomorphism
$(\ast)$ follows from  \cite[Cor. 3.2.,~]{P}.
Since $(\Z_2)^{n-1}$ is an abelian group,
$$\phi_{\Gamma} = \bigoplus_{i=1}^{n}(\phi_{\Gamma})_i,$$
where $(\phi_{\Gamma})_i:(\Z_2)^{n-1}\to\{\pm 1\}$ are irreducible representations of
$(\Z_2)^{n-1},$ for $i=1,2,...,n.$
From Lemma 1 and \cite[Theorem 8.2.2]{DH} the tangent bundle
$$\tau(M^n) = B(h)^{\ast}(\eta_n) = \bigoplus_{i=1}^{n} B(h_i)^{\ast}(\eta_1),$$
where $h_i = (\phi_{\Gamma})_{i}\circ h.$
This finishes the proof.

\pe
\newtheorem{cor}{Corollary}
\vskip 2mm
\section{$\operatorname{Spin}^{\C}$-structure}
%We keep notation from the previous sections.
It is well known (see \cite[Example 4.6 on page 4593]{MP}) that HW-manifolds do not have the $\operatorname{Spin}$-structure.
In this section we shall consider the question:
do HW-manifolds have the $\operatorname{Spin}^{\C}$-structure?
\vskip 1mm
\noindent
On the beginning let us recall some facts about the group $\operatorname{Spin}^{\C}$
defined in the Introduction.
We start with homomorphisms (\cite[page 25]{TF}):
\begin{itemize}
\item $i: \operatorname{Spin}(n)\to \operatorname{Spin}^{\C}(n)$ is the natural inclusion $i(g) = [g,1].$
\item $j:S^{1}\to \operatorname{Spin}^{\C}(n)$ is the natural inclusion, $j(z) = [1,z].$
\item $l:\operatorname{Spin}^{\C}(n)\to S^1$ is given by $l[g,z] = z^2.$
\item $p:\operatorname{Spin}^{\C}(n)\to \operatorname{SO}(n)\times S^1$ is given by $p([g,z]) = (\lambda_n(g),z^2).$
Hence $p = \lambda_n\times l.$
\end{itemize}
\vskip2mm
Since $S^1 = \operatorname{SO}(2),$ there is the natural map
$k:\operatorname{SO}(n)\times \operatorname{SO}(2)\to\operatorname{SO}(n+2).$
Then we can describe $\operatorname{Spin}^{\C}(n)$ as the pullback by~this~ map of the covering map
$$\begin{diagram}
\node{\operatorname{Spin}^\C(n)}\arrow{s,r}{p}\arrow{e}\node{\operatorname{Spin}(n+2)}\arrow{s,r}{\lambda_{n+2}}\\
\node{\operatorname{SO}(n)\times \operatorname{SO}(2)}\arrow{e,t}{k}\node{\operatorname{SO}(n+2)}
\end{diagram}$$
\vskip 5mm
\noindent
Let $W^n$ be an $n$-dimensional, compact oriented manifold and
let $\delta:W^n\to\operatorname{BSO}(n)$ be the classification map of its tangent bundle $TW^n.$
We now recall the definition of a $\operatorname{Spin}^{\C}$-structure (\cite[page 34]{K}, \cite[page 47]{TF}).
\newtheorem{defin}{Definition}
\begin{defin}\label{spinc}
A $\operatorname{Spin}^{\C}$-structure
on the manifold $W^n$ is a lift of $\delta$ to $\operatorname{BSpin}^{\C}(n),$ giving
a commutative diagram{\rm:}
$$\begin{diagram}
\node{}\node{\operatorname{BSpin}^\C(n)}\arrow{s,r}{B\left(\bar{\lambda}_n\right)}\\
\node{W^n}\arrow{ne}\arrow{e,t}{\delta}\node{\operatorname{BSO}(n).}
\end{diagram}$$
\end{defin}
\vskip 20mm
\begin{rem}\label{char}\mbox{   }{\rm
\begin{enumerate}
\item  (See \cite[Remark, d) on page 49]{TF}.) $W^n$ has the $\operatorname{Spin}^{\C}$-structure
if and only if the Stiefel-Whitney class $w_2\in H^2(W^n,\Z_2)$ is  the $\Z_{2} $-reduction of an
integral Stiefel-Whitney class $\tilde{w_2}\in H^2(W^n,\Z).$
\vskip 1mm
\noindent
\item  Let $K(\Z,2)$ and $K(\Z_2,2)$ be the Eilenberg-Maclane spaces.
From the homotopy theory $$H^2(W^n,\Z) = [W^n, BS^1] =[W^n, K(\Z,2)]$$ and
$H^2(W^n,\Z_2) = [W^n, K(\Z_{2},2)].$ Hence the above condition defines a commutative diagram
$$\begin{diagram}
\node{}\node{K(\Z,2)}\arrow{s}\\
\node{W^n}\arrow{ne,l}{\tilde{w}_2}\arrow{e,t}{w_2}\node{K(\Z_2,2)}
\end{diagram}$$
where the vertical arrow is induced by an epimorfizm $\Z\to\Z_2.$
\end{enumerate}}
\end{rem}
\vskip 3mm
\noindent
From previous sections (Lemma \ref{mainprop}) an oriented flat manifold $M = \operatorname{B}\Gamma,$ and
$\delta = \operatorname{B}(h)$ where $\Gamma = \pi_1(M)$ and $h:\Gamma\to \operatorname{SO}(n)$ is a holonomy
homomorphism. Let us recall (see\cite{TF}, \cite{MP1}, \cite{PS} and Remark~\ref{SW2} below) that an oriented manifold $M$ has a Spin-structure
if and only if there exists a homomorphism $e:\Gamma\to \operatorname{Spin}(n)$ such that
\begin{equation}\label{spin}
\lambda_n\circ e = h.
\end{equation}
Hence, a condition of existence of the $\operatorname{Spin}^{\C}$-structure on $M$ is very similar to
the above condition~(\ref{spin}).
\begin{theo}\label{main}
Let $M$ be an oriented flat manifold with $H^2(M,\R) = 0.$ $M$ has a $\operatorname{Spin}^{\C}$-structure
if and only if
there exists a homomorphism
$\epsilon:\Gamma\to \operatorname{Spin}^{\C}(n)$ such that
\begin{equation}\label{sp}
\bar{\lambda}_n\circ\epsilon = h.
\end{equation}
\end{theo}
{\bf Proof:}
Let us assume that there exists a homomorphism $\epsilon:\Gamma\to \operatorname{Spin}^{\C}(n)$
such that $\bar{\lambda}_{n}\epsilon = h.$
We claim that the conditions of Definition~\ref{spinc} are satisfied.
In fact, $B(\bar{\lambda}_{n})B(\epsilon) = B(h)$ up to homotopy.
To go the other way,
let us assume that $M = B\Gamma$
admits a $\operatorname{Spin}^{\C}$-structure.
We have a commutative diagram.
\vskip 10mm
\begin{center}
\newlength{\rsep}
\newlength{\csep}
\setlength{\rsep}{-3em}
\setlength{\csep}{3.5em}
\begin{tikzpicture}
%[every node/.style={execute at begin node=$, execute at end node=$}]
\path node (r11) at (\csep,0) {0} node (r12) at (3\csep,0) {$Z_2$} node (r13) at
(5\csep,0) {$\Gamma_0$} node (r14) at (7\csep,0) {$\Gamma$} node (r15) at (9\csep,0) {0};
\path node (r21) at (0\csep,\rsep) {0} node (r22) at (2\csep,\rsep) {$S^1$} node (r23)
at (4\csep,\rsep) {$\Gamma_1$} node (r24) at (6\csep,\rsep) {$\Gamma$} node (r25) at
(8\csep,\rsep) {0};
\path node (r31) at (\csep,2\rsep) {0} node (r32) at (3\csep,2\rsep) {$Z_2$} node
(r33) at (5\csep,2\rsep) {$\operatorname{Spin}(n)$} node (r34) at (7\csep,2\rsep) {$\operatorname{SO}(n)$} node (r35)
at (9\csep,2\rsep) {0};
\path node (r41) at (0\csep,3\rsep) {0} node (r42) at (2\csep,3\rsep) {$S^1$} node
(r43) at (4\csep,3\rsep) {$\operatorname{Spin}^{\C}(n)$} node (r44) at (6\csep,3\rsep) {$\operatorname{SO}(n)$} node
(r45) at (8\csep,3\rsep) {0};
\path[->] (r11) edge (r12) (r12) edge (r13) (r13) edge (r14) (r14) edge (r15);
\path[->] (r21) edge (r22) (r22) edge (r23) (r23) edge (r24) (r24) edge (r25);
\path[->] (r41) edge (r42) (r42) edge node[above] {j} (r43) (r43) edge node[above] {$\bar{\lambda}_n$}(r44) (r44) edge (r45);
\path[->] (r12) edge node[above] {r} (r22) (r13) edge (r23) (r14) edge[-,double] (r24);
\path[->] (r22) edge[-,double] (r42) (r23) edge (r43) (r24) edge (r44);
\path[->] (r32) edge node[above] {r} (r42) (r33) edge node[above] {i} (r43) (r34) edge[-,double] (r44);
\path node (e21) at (3\csep,\rsep) {\;} node (e22) at (5\csep,\rsep) {\;} node (e23)
at (7\csep,\rsep) {\;};
\path[->] (r12) edge[-,double] (e21) (e21) edge[-,double] (r32) (r13) edge[-] (e22)
(e22) edge (r33) (r14) edge[-] node[auto]{h} (e23) (e23) edge (r34);
\path node (e31) at (2\csep,2\rsep) {\;} node (e32) at (4\csep,2\rsep) {\;} node
(e33) at (6\csep,2\rsep) {\;};
\path[->] (r31) edge[-] (e31) (e31) edge (r32) (r32) edge[-] (e32) (e32) edge (r33)
(r33) edge[-] (e33) (e33) edge node[above]{$\lambda_n$}(r34) (r34) edge (r35);
\end{tikzpicture}
\end{center}
$$\text{Diagram 1}$$
\noindent
where $\Gamma_0$ is defined by the second Stiefel-Whitney class $w_2\in H^2(\Gamma,\Z_2)$ and
$\Gamma_1$ is defined by the element $r_{\ast}(w_2)\in H^2(\Gamma,S^1).$ Here $r:\Z_2\to S^1$
is a group monomorphism.
Let $h^2:H^2(\operatorname{SO}(n),K)\to H^2(\Gamma,K)$ be a homomorphism induced by the holonomy homomorphism $h,$
for $K = \Z_2, S^1.$
By definition, (see \cite[Chapter 23.6]{PM}) there exists an element\footnote{$H^{\ast}(\operatorname{SO}(n),\Z_2) = \Z_2[x_2,x_3,...,x_n].$}
$x_2\in H^2(\operatorname{SO}(n),\Z_2) = \Z_2$
such that $h^2(x_2) = w_2$ and
$h^2(r_{\ast}(x_2)) = r_{\ast}(h^2(x_2)) = r_{\ast}(w_2).$
Moreover, we have two infinite sequences of cohomology
which are induced by the following commutative diagram of groups:
\vskip 1mm
$$\begin{diagram}
\node{1}\arrow{e}\node{\Z}\arrow{s,=}\arrow{e,t}{2}\node{\Z}\arrow{s}\arrow{e}\node{\Z_2}\arrow{s,l}{r}\arrow{e}\node{1}\\
\node{1}\arrow{e}\node{\Z}\arrow{e}\node{\R}\arrow{e}\node{S^1}\arrow{e}\node{1}
\end{diagram}$$
\dgARROWLENGTH=0.5cm
$$\begin{diagram}
\node{\ldots}\arrow{e}\node{H^2(\Gamma,\Z)}\arrow{s,=}\arrow{e}\node{H^2(\Gamma,\Z)}\arrow{s}\arrow{e,t}{red}\node{H^2(\Gamma,\Z_2)}
\arrow{s,l}{r_*}\arrow{e}\node{H^3(\Gamma,\Z)}\arrow{s,=}\arrow{e}\node{\ldots}\\
\node{\ldots}\arrow{e}\node{H^2(\Gamma,\Z)}\arrow{e}\node{H^2(\Gamma,\R)}\arrow{e}\node{H^2(\Gamma,S^1)}\arrow{e}\node{H^3(\Gamma,\Z)}\arrow{e}
\node{\ldots}
\end{diagram}$$
From Remark~\ref{char} $\operatorname{red}(\tilde{w_2}) = w_2$ and since $H^2(\Gamma,\R) = 0, r_{\ast}(w_2) = 0.$ It follows that the row
$$0\to S^1\to\Gamma_1\to\Gamma\to 0$$
of the Diagram 1 splits. Hence there exists a homomorphism $\epsilon:\Gamma\to\operatorname{Spin}^{\C}(n)$ which satisfies~(\ref{sp}).
This proves the theorem.

\pe

\noindent
As an immediate corollary we have.
\begin{cor}
Let $M$ be an oriented flat manifold with the fundamental group $\Gamma$. If there exists a homomorphism
$\epsilon:\Gamma\to \operatorname{Spin}^{\C}(n)$ such that
\begin{equation}
\bar{\lambda}_n\circ\epsilon = h,
\end{equation}
then $M$ has a $\operatorname{Spin}^{\C}$-structure.
\vskip 5mm
\noindent
\end{cor}
\vskip 1mm
\begin{rem}\label{SW2}
Condition {\em (\ref{spin})} of the existence of the
$\operatorname{Spin}$-structure for oriented flat manifolds also follows from the proof of {\em Theorem~\ref{main}}.
\end{rem}
\vskip 1mm
\noindent
{\bf Question:} Is the assumption $H^2(M,\R)=0$ about the second cohomology group  necessary?
\vskip 1mm
\newtheorem{ex}{Example}
\begin{ex}\label{sissc}\mbox{   }
\begin{enumerate}
\item Because of the inclusion $i:\operatorname{Spin}(n)\to \operatorname{Spin}^{\C}(n)$ each $\operatorname{Spin}$-structure
on $M$ induces a $\operatorname{Spin}^{\C}$-structure.
\item If $M$ is any smooth compact manifold with an almost complex structure,
then $M$ has a canonical $\operatorname{Spin}^{\C}$-structure, see {\em \cite[page 27]{TF}}.
\end{enumerate}
\end{ex}
\begin{ex}\label{four}
Any oriented compact manifold of dimension up to four has a $\operatorname{Spin}^{\C}$-structure,
see {\em \cite[page 49]{G}}.
\end{ex}
From Example \ref{four} and \cite{RT} we have immediately.
\begin{cor}\label{nosyesc}
There exist three four dimensional flat manifolds without a {\rm Spin}-structure but
with a $\operatorname{Spin}^{\C}$-structure.
\end{cor}
In \cite[Example on page 50]{TF} is given a compact $5$-dimensional manifold
$Q$, without the $\operatorname{Spin}^{\C}$-structure. However the fundamental group $\pi_{1}(Q) = 1.$
There are also two other non-simply connected $5$-dimensional examples, see \cite[Eaxmples page 438]{KR}.
The first one is hypersurface in $\R P^2\times\R P^4$ defined by the equation $x_{0}y_{0}+x_{1}y_{1}+x_{2}y_{2} = 0$
where $[x_{0}:x_{1}:x_{2}]$ and $[y_{0}:y_{1}y_{2}:y_{3}:y_{4}]$ are homogeneous
coordinates in $\R P^2$ and $\R P^4$ respectively. The second is the Dold manifold
$$P(1,2) = \C P^2\times S^1/\sim,$$ where $\sim$ is an involution, which acts on $\C P^2$ by
complex conjugation and antipodally on $S^1.$
Our next result gives examples of $5$-dimensional flat manifolds without $\operatorname{Spin}^{\C}$-structure.
\begin{prop}\label{fived}
The {\rm HW}-manifolds $M_1$ and $M_2$ of dimension five do not have the $\operatorname{Spin}^{\C}$-structure.
\end{prop}
{\bf Proof:}
Since $H^2(M_{i},\R) = 0, i = 1,2,$ (\cite{DP}, \cite{RS}) we can apply a condition from Theorem~\ref{main}.
Let $\Gamma_{1} = \pi_1(M_1).$
It has the CARAT number 1-th 219.1.1, see \cite{carat}.\footnote{Here we use the name CARAT for tables of Bieberbach groups of dimension
$\leq 6,$ see \cite{carat}.}
It is generated by
$$\alpha_1 = ([1,1,1,-1,-1], (0,0,1/2,1/2,0)),
\alpha_2 = ([1,1,-1,-1,1],(0,1/2,0,0,0)),$$
$$ \alpha_3 = ([-1,1,1,-1,1],(0,0,0,0,1/2)),
\alpha_4 = ([1,-1,-1,1,1],(1/2,0,0,0,0))$$ and translations.
We assume that there exists a homomorphism $\epsilon:\Gamma_1\to \operatorname{Spin}^{\C}(5)$
such that $\bar{\lambda}_{n}\circ\epsilon = h.$
By definition
$$\alpha_2\alpha_3 = \alpha_3\alpha_2$$ and
$(\alpha_2\alpha_3)^2 = (\alpha_2)^2(\alpha_3)^2.$
Put $\epsilon(\alpha_i) = [a_i,z_i]\in \operatorname{Spin}^{\C}(5), a_i\in \operatorname{Spin}(5), z_i\in S^1, i = 1,2,3.$
Then
$$\epsilon\left((\alpha_2\alpha_3)^2\right) =\left [-1,z_{2}^{2}z_{5}^{2}\right]
= \epsilon((\alpha_2))^2\epsilon\left((\alpha_3)^2\right) =
\left[-1,z_{2}^{2}\right]\left[-1,z_{5}^{2}\right] = \left[1,z_{2}^{2}z_{5}^{2}\right]$$
and
$- z_{2}^{2}z_{5}^{2} = z_{2}^{2}z_{5}^{2}.$ We obtain a contradiction.

Now, let us consider the second five dimensional HW-group $\Gamma_2 = \pi_1(M_2)$ which
has a number 2-th. 219.1.1., (see \cite{carat}).
It is generated by
$$\beta_1 = (B_1,(1/2,1/2,0,0,0)), \beta_2 = (B_2, (0,1/2,1/2,0,0)),$$
$$\beta_3 = (B_3, (0,0,1/2,1/2,0))\hskip 2mm \text{and}\hskip 2mm \beta_4 = (B_4, (0,0,0,1/2,1/2)).$$
Put $\beta_5 = (\beta_1\beta_2\beta_3\beta_4)^{-1} = (B_5,(1/2,0,0,0,-1/2)).$
Assume that there exists a homomorphism $\epsilon:\Gamma_2\to \operatorname{Spin}^{\C}(5)$ which defines
the $\operatorname{Spin}^{\C}$-structure on $M_2.$ Let $\epsilon(\beta_i) = [a_i,z_i]\in \operatorname{Spin}^{\C}(5) = (\operatorname{Spin}(5)\times S^1)/\{1,-1\}.$
Let $t_i = (I,(0,...,0,\underbrace{1}_i,0,...,0)), i = 1,2,3,4,5.$
Since $\epsilon$ is a homomorphsim
\begin{equation}\label{five1}
\forall_{1\leq i\leq 5} \;\; \epsilon\left((\beta_i\beta_{i+2})^2\right) =
\left[a_{i}a_{i+2}a_{i}a_{i+2}, z_{i}^{2}z_{i+2}^{2}\right] = \left[-1,z_{i}^{2}z_{i+2}^{2}\right].
\end{equation}
Moreover, by an easy computation
\begin{equation}\label{five2}
\forall_{1\leq i\leq 5} \;\;(\beta_i)^2 = t_i, (\beta_{i}\beta_{i+2})^2 =
t_{i+1}t_{i+3}^{-1}\hskip 2mm\text{and}\hskip 2mm\epsilon(t_i) = \left[\pm 1,z_{i}^{2}\right].
\end{equation}
From (\ref{five1}), (\ref{five2}) and (\ref{ofour})
\begin{equation}\label{five3}
\left[-1,z_{1}^{2}z_{3}^{2}\right] = \left[1,z_{2}^{2}z_{4}^{2}\right] = \left[-1,z_{3}^{2}z_{5}^{2}\right] =
\left[1,z_{4}^{2}z_{1}^{2}\right] =\left[-1,z_{5}^{2}z_{2}^{2}\right] = \left[1,z_{1}^{2}z_{3}^{2}\right],
\end{equation}
which is impossible. Here indexes we read modulo $5.$
This finishes the proof.

\pe
\vskip 4mm
\begin{defin}\label{cyclic}
The {\rm HW}-manifold $M^n$ of dimension $n$ is cyclic if and only if
$\pi_1(M^n)$ is generated by the following elements {\em (see \cite[Lemma 1]{S1})}{\rm:}
$$\beta_i = (B_i, (0,0,0,...,0,\underbrace{1/2}_i,1/2,0,...,0)), 1\leq i\leq n-1,$$
$$\beta_n = (\beta_1\beta_2\dots\beta_{n-1})^{-1} = (B_n,(1/2,0,...,0,-1/2).$$
\end{defin}
We have
\begin{theo}\label{maintheorem}
Cyclic {\rm HW}-manifolds do not have the $\operatorname{Spin}^{\C}$-structure.
\end{theo}
{\bf Proof:} Since the above group $\Gamma_2$ satisfies our assumption
the proof is a generalization of arguments from Proposition~\ref{fived}.
Let $\Gamma$ be a fundamental group of the cyclic HW-manifold of dimension $\geq 5,$
with set of generators $\beta_i = (B_i,b_i), i = 1,2,.., n.$
Since $H^2(\Gamma,\R) = 0,$ (\cite{DP}) we can apply a condition from Theorem~\ref{main}.
Let us assume that there exist a homomorphism $\epsilon:\Gamma\to Spin^{\C}(n),$ which defines
the $\operatorname{Spin}^{\C}$-structure and
\begin{equation}\label{spinchom}
\epsilon(\beta_i) = [a_i,z_i], a_i\in \operatorname{Spin}(n), z_i\in S^1.
\end{equation}
From \cite{P} the maximal abelian subgroup $\Z^n$ of $\Gamma$ is exactly the commutator subgroup $\Gamma.$
Hence $\epsilon([\Gamma,\Gamma])\subset i\left(\operatorname{Spin}(n)\right)\subset \operatorname{Spin}^{\C}(n).$
Since $\forall_{i}\hskip 2mm \epsilon((\beta_i)^2) = [a_{i}^{2},z_{i}^{2}]$ and $(\beta_i)^2\in [\Gamma,\Gamma],$
$z_{i}^{2} = \pm 1,$ for $i = 1,2,...,n,$ it follows that
\begin{equation}\label{ofour}
\forall_{i}\hskip 2mm z_i \in\{\pm 1,\pm i\}.
\end{equation}
Let $t_i = (I,(0,...,0,\underbrace{1}_i,0,...,0)).$
From (\ref{spinchom})
\begin{equation}\label{one0}
\epsilon(t_i) = \epsilon\left((\beta_i)^2\right) = \left[\pm 1, z_{i}^{2}\right], i = 1,2,...,n
\end{equation}
and also
\begin{equation}\label{one1}
\forall_{1\leq i\leq n}\hskip 2mm \epsilon\left((\beta_i\beta_{i+2})^2\right) = \left[-1, z_{i}^{2}z_{i+2}^{2}\right].
\end{equation}
Moreover
\begin{equation}\label{one2}
\forall_{1\leq i\leq n}\hskip 2mm (\beta_i\beta_{i+2})^2 = t_{i+1}t_{i+3}^{-1}.
\end{equation}
From equations (\ref{one0}), (\ref{one1}) and (\ref{one2}) we have
$$\left[-1,z_{i}^{2}z_{i+2}^{2}\right] = \left[1,z_{i+1}^{2}z_{i+3}^{2}\right]$$
and
$$\forall_{1\leq i\leq n}\hskip 2mm z_{i}^{2}z_{i+2}^{2} = -z_{i+1}^{2}z_{i+3}^{2} = z_{i+2}^{2}z_{i+4}^{2}.$$
Since $n$ is odd $z_{i}^{2}z_{i+2}^{2} = - z_{i+n}^{2}z_{i+2+n}^{2} = - z_{i}^{2}z_{i+2}^{2},$ a contradiction,(cf. (\ref{five3}))~\footnote{The indexes should be read modulo $n.$}.
This finishes the proof.

\pe

\noindent\textbf{Acknowledgment}
\vskip2mm
We would like to thank J. Popko for his help in the proof of Theorem~\ref{main}, B. Putrycz
for discussion about the existence of the $\operatorname{Spin}^{\C}$-structures on HW-manifolds
and A. Weber for some useful comments.


\begin{thebibliography}{99}
\bibitem{AS} L. Auslander, R. H. Szczarba, Vector bundles over tori and noncompact solvmanifolds,
American J. Mathem. 97, (1975), pp. 260 - 281
\bibitem{Ch} Charlap L.S.: {\it Bieberbach Groups and Flat Manifolds.}
Springer-Verlag, 1986.
\bibitem{DR} K. Dekimpe, B. De Rock, The Anosov theorem for flat generalized Hantzsche - Wendt
manifolds, J. Geom. Phys. 52 (2004), No.2, 177 - 185
\bibitem{DP} K. Dekimpe, N. Petrosyan, Homology of Hantzsche-Wendt groups, Contemporary Mathematics, 501
Amer. Math. Soc. Providence, RI, (2009), 87 - 102
\bibitem{TF} T. Friedrich, {\it Dirac operators in Riemannian geometry}, Graduate Studies in Mathematics, Vol. 25,
American Mathematical Society, Providence, Rhode Island 2000
\bibitem{G} R. E. Gompf, $\operatorname{Spin}^{\C}$-structures and homotopy equivalences,
Geometry and Topology, 1, (1997), 41-45
\bibitem{DH} D. Husem\"oller, {\it Fibre bundles}, McGraw-Hill, New York 1966
%\bibitem{KT} F. Kamber, P. Tondeur, Flat bundles and characteristic classes of
%group-representations, Amer. J. Math. 89, (1967), 857 - 886
\bibitem{KR} T. P. Killingback, E. G. Rees, $\operatorname{Spin}^{\C}$-structures on manifolds, Class. Quantum Grav. 2 (1985), 433-438
\bibitem{K} R. C. Kirby, {\it The Topology of 4-Manifolds}, Springer LN 1374, New York 1989
\bibitem{PM} J. P. May, A concise course in Algebraic Topology, Chicago Lectures in Mathematics, University of Chicago Press,
Chicago 1999
\bibitem{MP} R. Miatello, R. Podest\'a, The spectrum of twisted Dirac operators on compact flat manifolds,
Trans. A.M.S., 358, Number 10, 2006, 4569 - 4603
\bibitem{MP1} R. Miatello, R. Podest\'a, Spin structures and spectra of $Z_2^k$-manifolds,
Math. Z., 247, 2004, 319-335
\bibitem{MR} R. Miatello, J. P. Rossetti, Isospectral Hantzsche-Wendt manifolds, J. Reine Angew. Math. 515 (1999),
1 - 23
\bibitem{carat}J. Opgenorth, W. Plesken, T. Schulz - CARAT - Crystallographic Algorithms and Tables -
http://wwwb.math.rwth-aachen.de/CARAT/
\bibitem{P}  B. Putrycz, Commutator Subgroups of  Hantzsche-Wendt Groups, J. Group Theory, 10 (2007), 401 - 409
\bibitem{PS} B. Putrycz, A. Szczepa\'nski: {Existence of spin structures on flat four - manifolds}
Adv. in Geometry, 10 (2), (2010), 323-332
\bibitem{RT} J. G. Ratcliffe, S. T. Tschantz, Spin and complex structures on flat gravitational instantons, Classical Quantum Gravity, 17 (2000), no. 1, 179 - 188
\bibitem{RS} J. P. Rossetti, A. Szczepa\'nski, Generalized Hantzsche-Wendt flat manifolds,
Revista Iberoam. Mat. {\bf 21}(3), 2005, 1053-1079
\bibitem{S0} A. Szczepa\'nski, Aspherical manifolds with the $\Q$-homology of a sphere,
Mathematika, 30, (1983), 291-294
\bibitem{S1} A. Szczepa\'nski, The euclidean representations of the Fibonacci groups,
{Q. J. Math.} {\bf 52} (2001), 385-389;
\bibitem{S2} A. Szczepa\'nski, Properties of generalized Hantzsche - Wendt groups, J. Group Theory
{\bf 12},(2009), 761-769
\bibitem{S3} A. Szczepa\'nski, {\it Geometry of Crystallographic Groups}, Algebra Ahd Discrete Mathematics, Vol. 4, World Scientific, 2012
\bibitem{ATV} A. T. Vasquez, Flat Riemannian manifolds, J. Diff. Geom. {\bf 4}, 1970, 367 - 382
\end{thebibliography}
\end{document}